\RequirePackage{ifpdf}
\ifpdf % We are running pdfTeX in pdf mode
\documentclass[pdftex]{sigma}
\else
\documentclass{sigma}
\fi

\usepackage{mathrsfs}
\usepackage{dsfont}

\newcommand{\set}[1]{\left\{#1\right\}}

\DeclareMathOperator{\spann}{\mathds{R}-span}
\DeclareMathOperator{\spannq}{\mathds{Q}-span}

 \DeclareMathOperator{\ind}{Ind}
\DeclareMathOperator{\Log}{log} \DeclareMathOperator{\Ad}{Ad}
\DeclareMathOperator{\ad}{ad}
 \DeclareMathOperator{\U}{\mathbf{U}}
\DeclareMathOperator{\V}{\mathbf{V}}

\DeclareMathOperator{\nongen}{\mathcal{V}}
\DeclareMathOperator{\gen}{\mathcal{W}}

\DeclareMathOperator{\Hi}{\mathscr{H}}

\DeclareMathOperator{\bbase}{\mathscr{B}}
\DeclareMathOperator{\kir}{\sans{Kir}}

\def\mult{\mbox{\usefont{T1}{cmfr}{m}{it}
  \selectfont m}}

\newfont{\Bfit}{cmbxti10 scaled\magstep1}

        \newfont{\Bfmit}{cmmib10 scaled\magstep1}

        \newfont{\Sans}{cmss10 scaled\magstep1}         \newcommand{\sans}[1]{\hbox{\Sans {#1}}}
        \newfont{\FRA}{eurm10 scaled\magstep1}
\def\dsr{{\mathds R}}
\def\dsc{{\mathds C}}
\def\dsn{{\mathds N}}
\def\dsq{{\mathds Q}}
\def\dsz{{\mathds Z}}

  \def\tore{{\mathbb I}}

\def\hb#1{\hbox{#1}}

\def\hb #1{\hbox{#1}}

\def\hb#1{\hbox{#1}}

\def\dim#1{\hb{dim}(#1)}
\def\Log#1{\rm{log}(#1)}

\def\L1#1{L^1(#1)}

\def\L#1#2{L^{#1}(#2)}

\def\lef({\left(}
\def\rig){\right)}

\begin{document}

\allowdisplaybreaks

\renewcommand{\PaperNumber}{021}

\FirstPageHeading

\ShortArticleName{Harmonic Analysis in One-Parameter Metabelian Nilmanifolds}

\ArticleName{Harmonic Analysis in One-Parameter\\ Metabelian Nilmanifolds}

\Author{Amira GHORBEL}

\AuthorNameForHeading{A.~Ghorbel}

\Address{Facult\'e des Sciences de Sfax,
D\'epartement de Math\'ematiques, \\
 Route de Soukra, B.P. 1171, 3000 Sfax, Tunisie}
\Email{\href{mailto:Amira.Ghorbel@fss.rnu.tn}{Amira.Ghorbel@fss.rnu.tn}}

\ArticleDates{Received September 02, 2010, in f\/inal form February 21, 2011;  Published online February 27, 2011}

\Abstract{Let $G$ be a connected, simply connected one-parameter
metabelian nilpotent Lie group, that means,  the corresponding Lie
algebra has a one-codimensional abelian sub\-algebra. In this article
we show that $G$ contains a discrete cocompact subgroup. Given a~discrete cocompact subgroup $\Gamma$ of $G$, we def\/ine the
quasi-regular representation $\tau = \ind_\Gamma^G 1$   of $G$. The
basic problem considered in this paper concerns the decomposition of
$\tau$ into irreducibles. We give an orbital description of the
spectrum, the multiplicity function and we construct an explicit
intertwining operator between $\tau$ and its desintegration  without
considering multiplicities. Finally, unlike the Moore inductive
algorithm for multiplicities on nilmanifolds, we carry out here a
direct computation to get the multiplicity formula.}


\Classification{22E27}

\section{Introduction}

Let $G$ be a  connected simply connected nilpotent Lie group having
a cocompact discrete subgroup~$\Gamma$, then the nilmanifold
$G/\Gamma$ has a unique invariant Borel probability measure $\nu$
and $G$ acts on $\mathbf{L}^2(G/\Gamma) = \mathbf{L}^2(G/\Gamma,
\nu)$ by the quasi-regular representation $\tau = \ind_\Gamma^G 1 $.
This means that
\begin{gather*}%\label{regular_representation}
\big(\tau(a) f\big)(g) = f(a^{-1} g),\qquad  f\in
\mathbf{L}^2(G/\Gamma),\qquad a, g \in G.
\end{gather*}
It is known \cite[p.~23]{Gelfand1} that the representation $\tau$
splits into a discrete direct sum of a countable number of
irreducible unitary representations, each $ \pi $ with f\/inite
multiplicity $\mult(\pi)  $. We write
\begin{gather}\label{disintegration-gelfand}
\tau \simeq   \sum _{\pi \in (G:\Gamma)} \mult(\pi)
 \pi.
\end{gather}
For any given nilpotent Lie group $G$ with discrete cocompact
subgroup $\Gamma$, there are two general problems to consider. The
f\/irst is to determine the spectrum $(G:\Gamma)$ and the multiplicity
function $\mult(\pi)$ of these representations. The second is to
construct an explicit  intertwining operator
  between $\tau$ and its decomposition into irreducibles.

  A necessary and suf\/f\/icient condition for
$\pi$ to occur in $\tau$ was obtained, by Calvin C.~Moore~\cite{Moore1}, in the special case in which $\Gamma$ is a~lattice
subgroup of $G$ (i.e., $\Gamma$ is a~discrete cocompact subgroup of
$G$ and $\log(\Gamma)$ is an~additive subgroup of the Lie algebra
$\frak g$ of $G$). Moreover, Moore determined an~algorithm which
expresses $\mult(\pi)$ in terms of multiplicities of certain
representations of a~subgroup of codimension one in~$G$. Later, R.~Howe and L.~Richardson \cite{Howe2,Richardson1}
independently obtained a closed formula for the multiplicities for
general $\Gamma$. In~\cite{Moore2}, C.~Moore and J.~Wolf determine
explicitly  which square integrable representations of $G$ occur in
the decomposition (\ref{disintegration-gelfand}) and give a method
for calculating the multiplicities. Using the Poisson summation and
Selberg trace formulas,  L.~Corwin and F.~Greenleaf gave a formula
for $\mult(\pi)$ that depended only on the coadjoint orbit in~$\frak
g^*$ corresponding to $\pi$ via Kirillov theory and the structure of
$\Gamma$~\cite{Cor9} (see also \cite{Cor4,Cor5, Cor7,Cor9,Cor10}).

Our main goals in this paper are twofold. The f\/irst is to give an
orbital description of the decomposition of $\tau$ into
irreducibles, in the case when $G$ is a connected, simply connected
one-parameter metabelian nilpotent Lie group. Such a decomposition
has two components, the spectrum $(G:\Gamma)$ and the multiplicity
function $\pi \mapsto \mult(\pi)$. We give orbital description
of both. The second main goal is to give an explicit intertwining
operator between $\tau$ and its desintegration.

 This
paper is organized as follows. In Section~\ref{section2}, we establish
notations and recall a few standard facts about representation
theory, rational structure and cocompact subgroups of a~connected
simply connected nilpotent Lie groups. Section~\ref{section3} is devoted to
present some results which will be used in the next sections. In
Section~\ref{section4} we prove, f\/irst, that a one-parameter metabelian
nilpotent Lie group admits a discrete uniform subgroup $ \Gamma$
(i.e., the homogeneous space $G/ \Gamma$ is compact). Next, for a~f\/ixed discrete uniform subgroup $\Gamma$ we prove the existence of a~one-codimensional abelian rational ideal $M$  (i.e., $M \cap \Gamma
$ is a~discrete uniform subgroup of~$M$). Furthermore, we give a
necessary and suf\/f\/icient condition for the uniqueness of~$M$. In
Section~\ref{section5}, we pick from a strong Malcev basis strongly based on
$\Gamma$, an orbital description of the spectrum $(G:\Gamma)$. We
obtain the following decomposition
\begin{gather*}%\label{decomposition-special}
\tau \simeq \rho =   \bigoplus_{l \in \Sigma}\ind_M^G\chi_l +
\bigoplus_{l\in \nongen} \chi_l,
\end{gather*}
 where   $\Sigma$ is a crosssection for
$\Gamma$-orbits in a certain subspace  $\gen \subset \frak g^*$ and
$\nongen\subset \frak g^*$ (more details in Theorem~\ref{des-operator-special}).

We describe an intertwining operator $\U$ of $\tau $ and $\rho$,
  def\/ined for all  $\xi\in \mathcal{C}(G/\Gamma)$ and $g\in G$ by
\[
 \U(\xi)(l)(g) =
    \int_{M(l)/M(l)\cap \Gamma} \xi(gm) \chi_l(m) d\dot{m},
\]
where $M(l)= M$ if  $l\in \Sigma$ or $ M(l) = G $ if $ l\in\nongen$.
This operator does not take into account the
 multiplicities of the decomposition of~$\tau$.
As a consequence, we give an orbital  description of the
multiplicity function.

\section{Notations and basic facts} \label{section2}

The purpose of this section is to establish notations which will be
used in the sequel and recall some basic def\/initions and results
which we shall freely use afterwards without mentioning them
explicitly.

 Let $G$  be a connected and simply
connected nilpotent Lie group with Lie algebra $\frak g$, then the
exponential map ${\rm exp}: \frak g \rightarrow G$ is a
dif\/feomorphism. Let $\log: G \rightarrow \frak g$
denote the inverse of~${\rm exp}$.

\subsection{Induced representation}\label{section2.1} %\label{induced-representation}

Starting from a closed
subgroup $H$ of a nilpotent Lie group $G$ and a unitary
representation $\sigma$ of~$H$ in a Hilbert space~$\Hi_\sigma$, let
us construct a unitary representation of~$G$. We realize the unitary
representation $\ind_H^G\sigma$ of~$G$, the representation induced
on~$G$ from the representation $\sigma$ of~$H$, by left translations
on the completion of the space of continuous functions~$F$ on~$G$
with values in~$\Hi_\sigma$ satisfying
\begin{gather}\label{covarience}
F(gh) = \sigma(h^{-1})(F(g)),\qquad  g\in G, \quad h\in H,
\end{gather}
and having a compact support modulo $H$, provided with the norm
\begin{gather*}%\label{norm}
\|F\| = \left(\nu_{G, H}\big(\|F\|^2\big)\right)^{\frac{1}{2}} =
\left(\int_{G/H} \|F(g)\|^2 d\nu_{G, H}(g)\right)^{\frac{1}{2}},
\end{gather*}
where
\begin{gather*}
\big(\big(\ind_H^G\sigma\big)(g) F\big)(x) = F\big(g^{-1} x\big),\qquad g, x\in G.
\end{gather*}

\subsection{The orbit theory}\label{section2.2}

Suppose $G$ is a nilpotent Lie group with Lie algebra $\frak g$. $G$~acts on $\frak g$ (respectively $\frak g^*$) by the adjoint
(respectively co-adjoint) action. For $l\in \frak g^*$, let
\[
\frak
g(l) = \{X\in \frak g:\; \langle l, [X, \frak g]\rangle = \{0\}\}
\]
be the stabilizer of $l$ in $\frak g$ which is actually the Lie
algebra of the Lie subgroup
\[
G(l) = \{g\in G:\ g\cdot l = l\}.
\]
So, it is clear that $\frak g(l)$ is the radical of the
skew-symmetric bilinear form $B_l$ on $\frak g$ def\/ined by
\begin{gather*}%\label{skew-symmetric bilinear}
B_l(X, Y) = \langle l, [X, Y]\rangle,\qquad X, Y \in \frak g.
\end{gather*}
A subspace $\frak b(l)$ of the Lie algebra $\frak g$ is called a
polarization for $l\in \frak g^*$ if it is a maximal dimensional
isotropic subalgebra with respect to $B_l$. We have the following
equality
\begin{gather}\label{dim-pilarisation}
\dim{\frak b(l)}=\dfrac{1}{2} \big(\dim{\frak g}+\dim{\frak
g(l)}\big).
\end{gather}

If $\frak b(l)$ is a polarization for $l\in \frak g^*$ let
$B(l)=\exp(\frak b(l))$ be the connected subgroup of $G$ with Lie
algebra $\frak b(l)$ and def\/ine a character of $B(l)$ by the formula
\begin{gather*}%\label{Defn-charcter}
\chi_l(\exp(X))= e^{2i \pi \langle l, X\rangle},\qquad \forall\, X\in \frak b(l).
\end{gather*}

Now, we recall the Kirillov orbital parameters. We denote by $\hat
G$ the unitary dual of $G$, i.e.\ the set of all equivalence classes
of irreducible unitary representations of $G$. We shall sometimes
identify the equivalence class $[\pi]$ with its representative $\pi$
and we denote the equivalence relation between two representations~$\pi_1$ and $\pi_2$ by $\pi_1 \simeq \pi_2$ or even by $\pi_1 =
\pi_2$.  The dual space $\hat G$ of $G$ is parameterized canonically
by the orbital space $\frak g^*/G$. More precisely, for~$l \in \frak
g^*$ we may f\/ind a real polarization $\frak b(l)$ for~$l$.
 Then the representation $\pi_l =
\ind_{B(l)}^G \chi_l$ is irreducible; its class is independent of
the choice of $\frak b(l)$; the Kirillov mapping
\begin{gather*}
\kir_G:  \ \ \frak g^* \rightarrow \hat G,\qquad  l \mapsto\pi_l
\end{gather*}
 is surjective and factors to a
bijection $\frak g^*/G \rightarrow \hat G$. Given $\pi \in \hat G$, we write $\Omega(\pi) \in \frak g^*/G$ to denote the inverse
image of $\pi$ under the Kirillov mapping $\kir_G$.

\subsection{Rational structures and uniform subgroups}\label{section2.3}

 In this section we present some results on discrete uniform
subgroups of nilpotent Lie groups.

%{\bf Rational structures.}
 Let $G$ be a nilpotent, connected
and simply connected real Lie group and let $\frak g$ be its Lie
algebra. We say that $\frak g$ (or $G$) has a \textit{rational
structure} if there is a Lie algebra $\frak g_\mathds{Q}$ over
$\mathds{Q}$ such that $\frak g \cong \frak g_\mathds{Q} \otimes
\mathds{R}$. It is clear that $\frak g$ has a rational structure if
and only if $\frak g$ has an $ \mathds{R}$-basis $(X_1,\dots,X_n)$
with rational structure constants.

%{\bf Uniform subgroups.}
 A discrete subgroup $\Gamma$ is
called \textit{uniform} in $G$ if the quotient space $G/\Gamma$ is
compact. The homogeneous space $G/\Gamma$ is called a
\textit{compact nilmanifold}. A proof of the next result can be
found in Theorem~7 of \cite{Malcev1} or in Theorem~2.12 of~\cite{Raghunathan}.

\begin{theorem}[the Malcev rationality criterion] \label{Malcev-cri}
Let $G$ be a simply connected nilpotent Lie
group, and let $\frak g$ be its Lie algebra. Then $G$ admits a
uniform subgroup $\Gamma$ if and only if $\frak g$ admits a basis
$(X_1,\dotsc,X_n)$ such that
\[
[X_i, X_j] =  \sum_{k=1}^n c_{ijk} X_k\qquad \forall\,  1\leq i, j\leq n,
\]
 where the constants $c_{ijk}$ are all
rational.
\end{theorem}

More precisely, if $G$ has a uniform subgroup $\Gamma$, then $\frak
g$  has a rational structure such that
\[
\frak g_\mathds{Q}= \frak g_{\mathds{Q}, \Gamma}= \spannq\set{\log(\Gamma)}.
\]
 Conversely, if
$\frak g$ has a rational structure given by some
$\mathds{Q}$-algebra $\frak g_\mathds{Q}\subset \frak g$, then~$G$
has a uniform subgroup $\Gamma$ such that $\log(\Gamma) \subset
\frak g_\mathds{Q}$ (see \cite{Cor1} or \cite{Malcev1}).

%{\bf Strong Malcev basis.}
\begin{definition}[\cite{Cor1}] Let $\frak g$ be a nilpotent Lie algebra and let
$\bbase=(X_1, \ldots, X_n)$ be a basis of $\frak g$.
 We say that $\bbase$ is a  strong  Malcev
basis for $\frak g$ if  $\frak g_i=\spann\set{X_1, \ldots, X_i}$ is
 an ideal of $\frak g$ for each $1\leq i\leq n$.
\end{definition}

  Let $
\Gamma$ be a discrete uniform subgroup of $G$. A strong Malcev basis
$(X_1,\dotsc,X_n)$ for $\frak g$ is said to be \textit{strongly
based on} $\Gamma$ if
\begin{gather*}%\label{strongly basis}
\Gamma= \exp(\mathds{Z}X_1)\dotsb \exp(\mathds{Z}X_n).
\end{gather*}
 Such a basis always exists (see \cite{Cor1,Matsushima1}).

\subsection{Rational subgroups}\label{2.4}

\begin{definition}[rational subgroup] Let $G$ be a connected
simply connected nilpotent Lie group with Lie algebra~$\frak g$. We
suppose that~$\frak g$ has rational structure given by $\frak
g_\dsq$.
\begin{itemize}\itemsep=0pt
  \item[$(1)$]
Let  $\frak h$ be an $\mathds{R}$-subspace of $\frak g$.  We say
that $\frak h$ is \textit{rational} if $\frak h = \spann\set{\frak
h_\mathds{Q}}$ where  $\frak h_\mathds{Q}= \frak h \cap \frak
g_\mathds{Q}$.
\item[$(2)$]
 A connected, closed subgroup
$H$ of $G$ is called \textit{rational} if its Lie algebra $\frak h$
is rational.
\end{itemize}
\end{definition}

\begin{remark}\label{ratio-intersection}
The $\dsr$-span and the intersection of rational subspaces are
rational \cite[Lemma 5.1.2]{Cor1}.
\end{remark}

\begin{definition}[subgroup with good $\Gamma$-heredity, \cite{Vinberg1}] Let $\Gamma$ be a discrete uniform subgroup in a
locally compact group $G$ and $H$ a closed subgroup in $G$. We say
that $H$ is a subgroup with good $\Gamma$-heredity if the
intersection $\Gamma\cap H$ is a discrete uniform subgroup of $H$.
\end{definition}

\begin{theorem}[\protect{\cite[Lemma A.5]{Cor9}}] Let $G$ be a connected, simply connected
nilpotent Lie group with Lie algebra  $\frak g$,  let $\Gamma$ be a
discrete uniform subgroup of $G$, and give $\frak g$ the rational
structure $\frak g_\dsq = \spannq\set{\log(\Gamma)}$.  Let $H$ be a
Lie subgroup of $G$. Then  the following statements are equivalent
\begin{itemize}\itemsep=0pt
  \item[$(1)$] $H$ is rational;
  \item[$(2)$] $H$ is a subgroup with good $\Gamma$-heredity;
  \item[$(3)$] The group $H$ is $\Gamma$-closed $($i.e., the set
  $H\Gamma$ is closed in $G)$.
\end{itemize}
\end{theorem}

A  proof of the next result can be found in Proposition 5.3.2 of~\cite{Cor1}.

\begin{proposition}\label{suite-adapte-rational}
Let $\Gamma$ be  discrete uniform subgroup in a nilpotent Lie group
$G$, and let $H_1\subsetneqq H_2\subsetneqq\cdots \subsetneqq H_k =G$
be  rational normal subgroups of $G$. Let $\frak h_1, \ldots, \frak
h_{k-1}, \frak h_k =\frak g$ be the corresponding Lie algebras. Then
there exists a strong Malcev basis $(X_1,\ldots, X_n)$ for $\frak g$
strongly based on $\Gamma$ and passing through $\frak h_1, \ldots,
\frak h_{k-1}$.
\end{proposition}

 Now let $G$ be a connected, simply connected nilpotent Lie
group with Lie algebra $\frak g$, and suppose that $\frak g$ has a
rational structure given by  the discrete uniform subgroup $\Gamma$.
A real linear functional $f \in \frak g^*$ is called rational ($f\in
\frak g_\mathds{Q}^*$,  $\frak
g_\mathds{Q}=\spannq\set{\Log{\Gamma}}$) if $\langle f, \frak
g_\mathds{Q}\rangle \subset \mathds{Q}$, or equivalently $\langle f,
\Log{\Gamma} \rangle \subset \mathds{Q}$.

\begin{proposition}[\cite{Cor9}, Theorem A.7] Let $G$ be a nilpotent Lie group with rational
structure and let $\frak g$ be its Lie algebra. If $l\in \frak g^*$
is rational, then its radical $\frak g(l)$ is rational.
\end{proposition}

\section{Preliminary  results}\label{section3}

\begin{definition} A functional $l\in \frak g^*$ is said in general
position or generic linear functional, if its coadjoint orbit has
maximum dimension.
\end{definition}
 The following proposition will be used in the sequel

\begin{proposition}\label{rational-functional}
Let $G = \exp(\frak g)$ be a nilpotent Lie group and $\Gamma$ a
discrete uniform subgroup of~$G$. Then there exist rational generic
linear functionals.
\end{proposition}

\begin{proof}
Let $\mathcal{O}$ be the set of elements in general position in
$\frak g^*$. We have $\mathcal{O}$ is a nonempty Zariski open set
in $\frak g^*$. Since $\frak g_\mathds{Q}^* $ is dense in $\frak
g^*$ then $\frak g_\mathds{Q}^*\cap \mathcal{O} \ne \varnothing$.
\end{proof}

%{\bf Fundamental domains for uniform subgroups.}
 Before stating the next result, we introduce some notations and
def\/initions. We f\/irst recall the concept of fundamental domain.

\begin{definition} Let $G$ be a topological group and let $H$ be a
subgroup of $G$. A fundamental domain for $H$ is a Borel subset
$\Omega$ of $G$ such that
\[
G= \bigsqcup_{h\in H} \Omega h
\] is
the disjoint union of the Borel subsets $ \Omega h$, $h\in H $.
\end{definition}

\begin{remark}
It is clear that a Borel subset  $\Omega$ of $G$ is a
fundamental domain for $H$ in $G$ if and only if  the natural map
$\Omega \rightarrow G/H$ def\/ined by $g\mapsto gH$ is
bijective.
\end{remark}

Let $\Gamma$ be a discrete uniform subgroup of a connected, simply
connected nilpotent Lie group~$G$, and let $ \bbase =(X_1, \ldots,
X_n)$ be a strong Malcev basis for the Lie algebra $\frak g$ of $G$
strongly based on~$\Gamma$. Def\/ine the mapping $\mathbf{E}_{\bbase}:
\mathds{R}^n\rightarrow G$ by
\begin{gather*}
\mathbf{E}_{\bbase}(T) = \exp(t_n X_n)\cdots \exp(t_1 X_1),
\end{gather*}
where $T= (t_1, \ldots, t_n)\in \mathds{R}^n$. It is well known that
 $\mathbf{E}_{\bbase}$ is a dif\/feomorphism~\cite{Cor1}. Let
\begin{gather*}%\label{tore}
\tore = [0, 1) = \set{t\in \mathds{R}:\ 0\leq t < 1}
\end{gather*}
and let
\begin{gather*}%\label{fundamental-domain}
\Omega= \mathbf{E}_{\bbase}(\tore^n).
\end{gather*}
Then $\Omega$ is a fundamental domain for $\Gamma$  in~$G$
\cite[Lemma~3.6]{Cor9}, and the mapping $\mathbf{E}_{\bbase}$ maps
the Lebesgue measure $dt$ on $\tore^n$ to the $G$-invariant
probability measure $\nu$ on $G/\Gamma$, that is, for~$\varphi $ in~$ \mathcal{C}(G/\Gamma)$, we have
\begin{gather*}%\label{measure-fundamental-domain}
\int_{G/\Gamma} \varphi( \dot{g}) d\nu( \dot{g}) = \int_{\tore^n}
\varphi( \mathbf{E}_{\bbase}(t)) dt.
\end{gather*}
Furthermore, for $\phi\in \mathcal{C}_c(G)$ we have
\begin{gather}\label{measure-fundamental-domain1}
\int_{G} \phi(g) d\nu(g) =  \sum_{s\in \dsz^n}\int_{\tore^n}
\phi\big( \mathbf{E}_{\bbase}(t) \mathbf{E}_{\bbase}(s)\big) dt.
\end{gather}

 The following
proposition will be used  in the sequel.

\begin{proposition}\label{decompositionmesure}
Let $G = \exp(\frak g)$ be a connected, simply connected nilpotent
Lie group and let~$\Gamma$ be a discrete uniform subgroup of $G$.
Let  $\frak m$ be a rational ideal of $\frak g$ of dimension $k$.
Let $M= \exp(\frak m)$ and let $(X_1, \ldots, X_n)$ be a strong
Malcev basis of $\frak g$ strongly based on $\Gamma$ passing through
$\frak m $.  For every $\xi\in \mathcal{C}_c(G/M)$, we have
\begin{gather*}%\label{100}
\int_{G/M}  \xi(g)   dg   =
 \sum_{s\in \mathds{Z}^{n-k}}
\int_{\tore^{n-k}}  \xi\Big(\exp(x_{n-k} \overline{X}_n)\cdots
\exp(x_1 \overline{X}_{k+1})\\
\hphantom{\int_{G/M}  \xi(g)   dg   =}{}\times
\exp(s_{n-k} \overline{X}_n)\cdots
\exp(s_1 \overline{X}_{k+1})\Big)
   dx_n\cdots dx_{k+1},
   \end{gather*}
where $\overline{X}_i$ are the image of $X_i$ under the canonical
projection $p: \frak g\rightarrow \frak g/\frak m$.
\end{proposition}

 \begin{proof}
Since $M$ is rational, it follows  from  Lemma 5.1.4 of \cite{Cor1},
that $P(\Gamma)$ is a uniform subgroup of $P(G)$, where $P: G
\rightarrow G/M$ is the canonical projection. On the other hand,
$(\overline{X}_{k+1}, \ldots, \overline{X}_{n})$ is a strong Malcev
basis of $\frak g/\frak m$ strongly based on $P(\Gamma)$. We
conclude by applying (\ref{measure-fundamental-domain1}).
 \end{proof}

\section{On the rational structure  of
one-parameter\\ metabelian nilpotent Lie groups}\label{section4}

\begin{definition}[one-parameter metabelian nilpotent Lie algebra]
A nonabelian nilpotent Lie algebra  $\frak g$ is said to be one
parameter metabelian nilpotent Lie algebra, if it admits  a
co-dimensional one abelian ideal
 in $\frak g$. \end{definition}

First, we give an important example of one-parameter metabelian
nilpotent Lie groups.

 \begin{example}[the generic f\/iliform nilpotent Lie group]
  Let $G$ be the generic f\/iliform
nilpotent Lie group of dimension $n$ with Lie algebra $\frak g$,
where
\begin{gather*}
 \frak g = \spann\set{X_1,
\ldots, X_n}
\end{gather*}
The Lie brackets  given by
\begin{gather*}
[X_n, X_i] = X_{i-1},  \qquad i=2,\ldots,  n-1,
\end{gather*}
and the nondef\/ined brackets being equal to zero or obtained by
antisymmetry. It is clear that~$\frak g$ is a one-parameter
metabelian nilpotent Lie algebra.
\end{example}

\begin{remark} Let $\frak g$ be a one-parameter metabelian
 nilpotent Lie algebra. Any
codimension one abelian ideal $\frak m\subset \frak g$ is a common
polarization for all functionals $l\in \frak g^*$ in general
position (i.e.,  $l\vert_{[\frak g, \frak g]}\ne 0$).
\end{remark}

\begin{definition}[one-parameter metabelian nilmanifold]\sloppy
A factor space of a one-parameter met\-abe\-lian nilpotent Lie group
over a discrete uniform subgroup is called a one-parameter
metabelian nilmanifold.
\end{definition}

\begin{proposition}\label{decomposition-speciale} Let $\frak g$ be a
one-parameter metabelian nilpotent Lie algebra. Then we have the
following decomposition
\begin{gather*}
\frak g = \mathds{R} X \oplus \oplus_{i=1}^p \mathcal{L}_{n_i}
\oplus \frak a
\end{gather*}
such that for all $i=1,\ldots, p$, the subalgebra $\mathds{R} X
\oplus  \mathcal{L}_{n_i} $ is the generic filiform nilpotent Lie
algebra of dimension $n_i+1$ and $\frak a \subset \frak z(\frak g)$.
 \end{proposition}

 \begin{proof} Let $\frak I$ be a one-codimensional abelian ideal of
 $\frak  g$ and let $X\in \frak g$ such that
\[
\frak g= \frak I\oplus \spann\set{X}.
\]
 Let $\mathrm{ad} X\vert_{\frak I}$ be the
 restriction of $\mathrm{ad}
 X$ to $\frak I$;
\[
\mathrm{ad}\, X\vert_{\frak I}: \quad \frak I\rightarrow \frak I,\qquad
 Y\mapsto [X, Y].
 \]
 Remark that $\mathrm{ad}\, X$ acts as a nilpotent linear
 transformation on~$\frak I$.
 By the Jordan normal form theorem, the matrix of
 $\mathrm{ad} X\vert_{\frak I}$ is similar to a matrix in real Jordan canonical form. Then
 there exist $\mathscr{B}_{\frak I}=(e_1, \ldots, e_{n-1})$  a
 basis of $\frak I$, $J_{0, n_i}$, $1\leq i\leq s$  elementary Jordan
 blocks of order $n_i$ such that
\[
\mathrm{Mat}(\mathrm{ad}\, X\vert_{\frak
 f}, \mathscr{B}_{\frak I}) =
  \mathrm{diag}[J_{0, n_1}, \ldots, J_{0, n_s}].
\]
  This completes the proof.
 \end{proof}

\begin{proposition}
Let $G$ be a one-parameter metabelian nilpotent Lie group with Lie
algebra~$\frak g$.
 Then $G$ admits a discrete uniform subgroup.
\end{proposition}

\begin{proof} This follows at once from
 the  Malcev  rationality criterion
 (Theorem~\ref{Malcev-cri}) and Proposition~\ref{decomposition-speciale}.
 \end{proof}

\begin{proposition}\label{maryem10}
Let $G$ be a one-parameter metabelian nilpotent Lie group of
dimension $n$.
 Let~$\Gamma $ be a discrete uniform subgroup of $G$. Then $G$ admits a
rational
 abelian ideal of codimension one.
\end{proposition}

Before proving the proposition, we need the following theorem.
\begin{theorem}[\protect{\cite[p.~17]{Bernat}}]
If $l\in \frak g^*$ is in general position then $\frak g(l)$ is
abelian.
\end{theorem}

\begin{proof}[Proof of Proposition~\ref{maryem10}]
Let $(e_1, \ldots, e_{n})$ be a strong Malcev basis for $\frak g$
strongly based on~$\Gamma$. Let $l$ be a rational  linear functional
in general position (see Proposition \ref{rational-functional}).
 The stabilizer of $l$ is a~rational abelian  subalgebra. On the other hand, applying equality~(\ref{dim-pilarisation}) we obtain
\begin{gather}\label{dim-stab}
\dim{\frak g(l)}=
 n-2.
 \end{gather}
  Let $(e_{i_1}, e_{i_2})$  be a basis of
$\frak g$ modulo $\frak g(l)$. If $ [ e_{i_1}, \frak g(l)] = \lbrace
0\rbrace$, then $\frak g(l)\oplus \spann\lbrace e_{i_1}\rbrace$  is
a rational abelian ideal of $\frak g$ of codimension one. If $ [
e_{i_1}, \frak g(l)] \ne \lbrace 0 \rbrace$, let
\[
\frak a = \lbrace X \in \spann\lbrace e_{i_1}, e_{i_2}\rbrace:\
[ X, \frak g(l)] =\set{0}\rbrace.
\] In this case, we have $\dim
{\frak a} =1$. Moreover, it is clear that
\[
 \frak a = \lbrace X \in \frak g:\
[ X, \frak g(l)] =\set{0}\rbrace\cap \spann\lbrace e_{i_1},
e_{i_2}\rbrace.
\]
 By Proposition~5 of~\cite{Ghorbel1}, the space
\[
\lbrace X \in \frak g:\
[ X, \frak g(l)] =\set{0}\rbrace
\] is rational. Consequently, $\frak
a$ is rational subspace in~$\frak g$. Thus $\frak a \oplus \frak
g(l)$ is a rational abelian ideal of~$\frak g$ of codimension one.
\end{proof}

The next proposition is a  generalization of  Proposition~3.1 of~\cite{Hamrouni1}, in which we give a necessary and suf\/f\/icient
condition for  uniqueness of the one-codimensional abelian normal
subgroup.

\begin{proposition}
 Let $G$ be a one-parameter metabelian nilpotent Lie group. Then $G$
admits a~unique one-codimensional
 abelian normal subgroup  if and only if $G$ is not of the form
 $H_3\times \dsr^k$, where $k\in \dsn$ and $H_3$ is the
 Heisenberg group  of dimension~$3$.
\end{proposition}

\begin{proof} The necessity of this condition is evident. We prove
the suf\/f\/iciency. Let $\frak m_1$, $\frak m_2$ be two abelian
one-codimensional subalgebras of~$\frak g$. Let $l$ be a linear
functional in general position. Since every polarization for $l$
contains~$\frak g(l)$ then there exists~$X\in \frak g$ (see~(\ref{dim-stab})) such that
\begin{gather*}
\frak m_1= \frak g(l)\oplus \spann\set{X}.
\end{gather*}
Let $Y\in \frak g$ such that
\begin{gather*}
\frak g=\frak m_1\oplus \spann\set{Y}.
\end{gather*}
On the other hand, there exist $\alpha\in \dsr^*$ and $u\in \frak
m_1$ such that
\begin{gather*}
\frak m_2= \frak g(l)\oplus \spann\set{\alpha Y+u}.
\end{gather*}
Since $\frak m_1$ and $\frak m_2$ are abelian then
\begin{gather*}
[Y, \frak g(l)]=\set{0}.
\end{gather*}
Consequently, $[X,\! Y]$ is  nonzero bracket. Let $\frak a$ be an
ideal of~$\frak g$ such that $\frak g(l)\! =\! \frak a\oplus
\spann\set{[X,\!Y]}$. Then we have
 $\frak g= \frak h_3\oplus \frak a$
 where $\frak h_3=\spann\set{ X,
Y , [X,Y]}$ is the three dimensional Heisenberg algebra.
\end{proof}

\section{Construction of intertwining operators}\label{section5}

Let $G$ be a one-parameter metabelian nilpotent  Lie group  of
dimension $n$ with Lie algebra~$\frak g$.
 Let $\Gamma$ be a discrete uniform
subgroup of~$G$. Let $M = \exp(\frak m)$, where $\frak m$ is  a
rational abelian ideal in~$\frak g$ of codimension one.
 Let  $\bbase = (X_1, \ldots, X_n)$ be a
strong Malcev basis of $\frak g$ strongly based on~$\Gamma$ passing
through $[\frak g, \frak g]$ and $\frak m$. We put
\begin{gather}\label{alg-derivee}
[\frak g, \frak g] = \spann\lbrace X_1,\ldots, X_p\rbrace
\end{gather}
and
\[
\frak m
=\spann\set{X_1, \ldots, X_{n-1}}.
\]
  Let
\[
\nongen = \mathds{Z} X_{p+1}^*+\cdots+\mathds{Z}X_{n}^*
\]
and
\[
\gen= \{ l \in
\mathds{Z}X_1^*+\cdots +\mathds{Z}X_{n-1}^*:\ l\vert_{[\frak g,
\frak g]}\ne 0\}.
\]
 In the following, for $l \in \frak m^* \subset \frak g^*$ and $g
\in G$ we denote
\[
 \mathrm{Ad}_0^*g l =(\mathrm{Ad}^* g l)\vert_{\frak m}.
 \]

\begin{lemma}
 The subset $\gen$ is $  \mathrm{Ad}_0^* \Gamma$ invariant.
\end{lemma}
\begin{proof}
Let $\exp(\gamma) \in \Gamma$ and $l \in \gen$. Let $i=1, \ldots,
n-1$. By def\/inition of the coadjoint representation we have
\[
\langle \exp(-\gamma)\cdot l, X_i\rangle = \langle l, e^{\ad
\gamma}(X_i)\rangle.
\]
 On the other hand, since $M$ is normal in
$G$,  we have  that
\[
\exp\big(e^{\ad
\gamma}(X)\big) = \exp(\gamma) \exp{X_i} \exp(-\gamma) \in \Gamma \cap
M.
\] Since $M$ is abelian and
\[
\Gamma \cap M=\exp(\dsz X_1)\cdots
\exp(\dsz X_{n-1})
\] then
\[
\log(\Gamma) \cap \frak m=\dsz X_1\oplus\cdots\oplus \dsz
X_{n-1}.
\] It follows that  $\langle l, e^{\ad \gamma}(X_i)\rangle
\in \dsz$. Then
\begin{gather}\label{cond1}
\exp(-\gamma)\cdot l\in
\dsz X_1^*\oplus\cdots\oplus \dsz X_{n-1}^*.
\end{gather}
It remains to prove that
\begin{gather}\label{cond2}
(\exp(-\gamma)\cdot l)\vert_{[\frak g, \frak g]}\ne 0.
\end{gather}
We have
\[
\frak g(\exp(-\gamma)\cdot l)=\Ad\exp(-\gamma)(\frak g(l)).
\]
Since $\frak g(l)\ne \frak g$, then $\frak g(\exp(-\gamma)\cdot
l)\ne \frak g$ and therefore (\ref{cond2}) holds. Consequently, from
(\ref{cond1}) and~(\ref{cond2}) we have $\exp(-\gamma)\cdot l \in
\gen$.
\end{proof}

Let $\mathcal{C}(G/\Gamma)$ be the space of all complex valued
continuous functions $\xi$ on $G$ satisfying
\begin{gather}\label{periode}
\xi(g\gamma)= \xi(g), \end{gather}
 for any $g$ in $G$ and $\gamma$
in $\Gamma$.

\begin{lemma}\label{periode2} Let $l\in \nongen$, we have
\begin{itemize}\itemsep=0pt
  \item[$(1)$] $\chi_l\vert_{\Gamma}=1$.
  \item[$(2)$] For $\xi \in \mathcal{C}(G/\Gamma)$, and $g\in G$, the function
\[
G/\Gamma\rightarrow \dsc,\qquad m\Gamma\mapsto \xi(gm)
\chi_l(m)
\]
is well defined.
\end{itemize}
\end{lemma}

\begin{proof} $(1)$ Let $\gamma\in \Gamma$, then there exist
$t_1, \ldots, t_n\in \dsz$ such that
\[
\gamma=\exp(t_1X_1)\cdots \exp(t_nX_n).
\]
By the Baker--Campbell--Hausdorf\/f formula we have
\[
\log(\gamma)\equiv t_{p+1}X_{p+1}+\cdots+t_n X_n\pmod{[\frak g,
\frak g]}.
\]
 Since $l\vert_{[\frak g, \frak g]}=0$ then $\chi_l(\gamma)=1$.

$(2)$ Let $\gamma\in \Gamma$. We have
\begin{gather*}
 \xi(gm\gamma)
\chi_l(m\gamma)
   = \xi(gm)
\chi_l(m\gamma) \overset{\text{by (\ref{periode})}}{=}    \xi(gm) \chi_l(m) \chi_l(\gamma)
 = \xi(gm) \chi_l(m).
\end{gather*}
 This completes the proof of the lemma.
\end{proof}

Let $\Sigma$ be a crosssection for $\Gamma$-orbits in $\gen$. Let
\begin{gather*}
 \rho =   \bigoplus_{l \in \Sigma}\ind_M^G\chi_l +
\bigoplus_{l\in \nongen} \chi_l.
\end{gather*}

We are now in a position  to formulate the following

\begin{theorem}\label{des-operator-special}
The operator  $\U$ defined for all  $\xi\in \mathcal{C}(G/\Gamma)$
and $g\in G$ by
 \begin{gather*}%\label{operator-special}
 \U(\xi)(l)(g) = \begin{cases}
\displaystyle     \int_{M/M\cap \Gamma} \xi(gm) \chi_l(m) d\dot{m} & \mbox{if }  l\in \Sigma,\vspace{1mm}\\
\displaystyle     \int_{G/ \Gamma} \xi(gm) \chi_l(m) d\dot{m} & \mbox{if }  l\in \nongen
    \end{cases}
\end{gather*}
 is an isometric linear operator having value in the Hilbert space
 $\Hi_\rho$ of $\rho$  and can be extended on
  $\mathbf{L}^2(G/\Gamma)$ to an
intertwining operator of $\tau$ and $\rho$.
\end{theorem}

\begin{proof} Clearly
for $ l \in \Sigma\cup \nongen$, the function $ \U(\xi)(l)$
satisf\/ies the covariance relation (\ref{covarience}). First, we
establish that $\U$ is well def\/ined and isometric. Let $\xi\in
\mathcal{C}(G/\Gamma)$. Then
\begin{gather*}
  \|\U(\xi)\|^2  = \sum_{l \in \Sigma}\|\U(\xi)(l)\|^2_{L^2(G/M, l)} + \sum_{l\in
  \nongen}
    \|\U(\xi)(l)\|^2.
\end{gather*}
We proceed to calculate the f\/irst sum. For $x=(x_1,\ldots, x_{n-1})
\in \mathds{R}^{n-1}$ and $ t\in \mathds{R}$, let
\begin{gather*}
\delta(t, x) = \mathbf{E}_{\bbase}((x, t))= \exp(t X_n) \exp(x_{n-1} X_{n-1})\cdots \exp(x_1
X_1),
\\
\sum_{l\in \Sigma} \|\U(\xi)(l)\|^2_{L^2(G/M, l)}
   = \sum_{l\in \Sigma} \int_{G/M} | \U(\xi)(l)(g)|^2     d\dot{g }\\
 \hphantom{\sum_{l\in \Sigma} \|\U(\xi)(l)\|^2_{L^2(G/M, l)}} {}
  =  \sum_{l\in \Sigma} \int_{G/M} \bigg| \int_{M/M\cap \Gamma}
\xi(gm) \chi_l(m) d\dot{m}\bigg|^2     d\dot{g}\\
\hphantom{\sum_{l\in \Sigma} \|\U(\xi)(l)\|^2_{L^2(G/M, l)}} {} = \sum_{l\in \Sigma} \sum_{s\in \mathds{Z}} \int_{\tore} \bigg|
\int_{M/M\cap \Gamma} \xi(\delta(t, 0)\delta(s, 0)  m) \chi_l (m)
 d\dot{m} \bigg|^2  dt \\
\hphantom{\sum_{l\in \Sigma} \|\U(\xi)(l)\|^2_{L^2(G/M, l)}} {} =
\sum_{l\in \Sigma} \sum_{s\in \mathds{Z}} \int_{\tore} \bigg|
\int_{M/M\cap \Gamma} \xi(\delta(t, 0)  m) \chi_l ( \delta(s,
0)^{-1} m \delta(s, 0))
 d\dot{m} \bigg|^2  dt\\
\hphantom{\sum_{l\in \Sigma} \|\U(\xi)(l)\|^2_{L^2(G/M, l)}} {}  =  \sum_{l\in \Sigma} \sum_{s\in
\mathds{Z}} \int_{\tore}  \bigg| \int_{M/M\cap \Gamma} \xi(\delta(t, 0)
m) \chi_{\Ad^*(\delta(s, 0))l} (m)
 d\dot{m} \bigg|^2  dt .
\end{gather*}
As the mapping $\Sigma\times \mathds{Z}\rightarrow \gen$, $(l, s) \mapsto \Ad^*(\delta(s, 0))l$
 is bijective, then
 \begin{gather*}
 \sum_{l\in \Sigma} \|\U(\xi)(l)\|^2_{L^2(G/M, l)}
   =  \sum_{l\in \gen}  \int_{\tore}  \bigg| \int_{M/M\cap
\Gamma} \xi(\delta(t, 0) m) \chi_l (m)
 d\dot{m} \bigg|^2  dt  \\
 \hphantom{\sum_{l\in \Sigma} \|\U(\xi)(l)\|^2_{L^2(G/M, l)}}{}
  =  \sum_{l\in \gen}  \int_{\tore}
 \bigg|
 \int_{\tore^{n-1}} \xi(\delta(t, x))
  \chi_l (\delta(0, x))
 dx  \bigg|^2  dt.
\end{gather*}
Next, we compute
\begin{gather*}
\sum_{l\in
  \nongen}
    \|\U(\xi)(l)\|^2  =
\sum_{l\in \nongen}
    | \U(\xi)(l)(e) |^2  = \sum_{l\in \nongen}\bigg|\int_{G/\Gamma}
    \xi(g) \chi_l(g) d\dot{g} \bigg|^2\\
\hphantom{\sum_{l\in   \nongen}     \|\U(\xi)(l)\|^2}
 = \sum_{l\in \nongen}\bigg|\int_{\tore^{n}}
    \xi(\delta(t, x)) \chi_l(\delta(t, x)
 ) dt dx \bigg|^2\\
\hphantom{\sum_{l\in   \nongen} \|\U(\xi)(l)\|^2}
 = \sum_{l_1\in \nongen_1} \sum_{l_2\in \nongen_2} \bigg|\int_{\tore^{n}}
    \xi(\delta(t, x)) \chi_{l_1}(\delta(0, x))
     \chi_{l_2}(\delta(t, 0) ) dt dx \bigg|^2\\
\hphantom{\sum_{l\in   \nongen} \|\U(\xi)(l)\|^2}
 (\text{where} \ \nongen_1 = \mathds{Z} X_{p+1}^*+\cdots+\mathds{Z}X_{n-1}^*,
     \nongen_2 =\mathds{Z}X_{n}^* \text{ and } l = l_1+l_2
      \in \nongen_1\oplus \nongen_2) \\
\hphantom{\sum_{l\in   \nongen} \|\U(\xi)(l)\|^2}
 = \sum_{l_1\in \nongen_1} \int_{\tore} \bigg|\int_{\tore^{n-1}}
    \xi(\delta(t, x))
 \chi_{l_1}(\delta(0, x))  dx \bigg|^2 dt,
\end{gather*}
where we have applied Parseval's equality with respect to the
variable $t$. Summarizing, we have
\begin{gather*}
  \|\U(\xi)\|^2    =  \sum_{l\in \gen}  \int_{\tore}
 \bigg|
 \int_{\tore^{n-1}} \xi(\delta(t, x))
  \chi_l (\delta(0, x))
 dx  \bigg|^2  dt \\
\hphantom{\|\U(\xi)\|^2    =}{}
+  \sum_{l\in \nongen_1} \int_{\tore}
\bigg|\int_{\tore^{n-1}}
    \xi(\delta(t, x))\chi_{l}(\delta(0, x))
dx \bigg|^2 dt.
\end{gather*}

It is clear that $\gen$ and $\nongen_1$ are disjoint. In fact,
suppose that $\gen\cap\nongen_1\ne \varnothing$. Let $l\in
\gen\cap\nongen_1$. The condition $l\in \nongen_1$ implies that
$l\vert_{[\frak g, \frak g]}= 0$ since $[\frak g, \frak g]=\dsr
X_1\oplus\cdots\oplus \dsr X_p$ (see~(\ref{alg-derivee})). This
contradicts the fact that $l\vert_{[\frak g, \frak g]}\ne 0$ because
$l\in \gen$. Now, let
\begin{gather}\label{disjoint-reunion}
\frak m_{\mathds{Z}}^* =\gen\sqcup \nongen_1.
\end{gather}
 It is
clear that
\[
\frak m_{\mathds{Z}}^* = \mathds{Z}X_1^*+\cdots
+\mathds{Z}X_{n-1}^*.
\]
 Consequently, we obtain
\begin{gather*} \|\U(\xi)\|^2
= \int_{\tore} \sum_{l\in \frak
m_{\mathds{Z}}^*} \left|\int_{\tore^{n-1}}
    \xi(\delta(t, x)) \chi_{l}(\delta(0, x))
dx \right|^2 dt   = \int_{\tore} \int_{\tore^{n-1}}
    |\xi(\delta(t, x))
 |^2 dx  dt \\
 \hphantom{\|\U(\xi)\|^2 }{}
 \overset{(\text{by Parseval's equality})}{=}  \int_{\tore^{n}}
    |\xi(\delta(t, x))
 |^2 dx  dt
   = \|\xi\|_{\mathscr{H}_\tau}^2.
\end{gather*}
The operator $\U$ being now isometric, it can be extended to an
isometry of $\mathbf{L}^2(G/\Gamma)$. It remains to verify that $\U$
is an intertwining operator for $\tau$ and $\rho$. It suf\/f\/ices then
to prove that  $\U\circ \tau(g)(\xi) = \rho(g) \circ \U(\xi)$ for
every  $g$ in $G$ and $\xi$ in $\mathcal{C}(G/\Gamma)$. Let $l\in
\Sigma, a\in G$. We compute
\[
\U\circ \tau(g)(\xi)(l)(a)=\int_{M/M\cap \Gamma} \tau(g)(\xi)(am)
\chi_l(m) d\dot{m}= \int_{M/M\cap \Gamma} \xi(g^{-1}am) \chi_l(m)
d\dot{m}.
\]
 On the other hand
\begin{gather*}
\rho(g) \circ \U(\xi)(l)(a)=\ind_M^G\chi_l(g)(\U(\xi)(l))(a)=
\U(\xi)(l)\big(g^{-1}a\big)\\
\hphantom{\rho(g) \circ \U(\xi)(l)(a)}{}
= \int_{M/M\cap \Gamma} \xi(g^{-1}am) \chi_l(m)
d\dot{m}.
\end{gather*}
Thus
\[
\U\circ \tau(g)(\xi)(l)(a)=\rho(g) \circ \U(\xi)(l)(a).
\]
Similarly, we prove that the same equality holds if $l\in \nongen$.
 Consequently,
 the following
diagram:
 \begin{gather}\label{intertwining}
\begin{CD} \mathscr{H}_\tau@>{\tau(g)}>> \mathscr{H}_\tau\\
@V{\U}VV  @V{\U}VV\\
\mathscr{H}_\rho @> {\rho(g)}>> \mathscr{H}_\rho
\end{CD}
\end{gather}
is commutative.
\end{proof}

Next, we show that $\U$ is an invertible operator. For this, let
\[
 \Hi_\rho^c =  \bigoplus_{l\in \Sigma} \mathcal{C}_c(G/M,
\chi_l) \oplus \bigoplus_{l\in \nongen} \mathds{C} \chi_l \subset
\mathscr{H}_\rho,
\]
 where $\mathcal{C}_c(G/M, \chi_l)$ is the space
of complex valued continuous functions $\xi$ on $G$ satisfying
\begin{gather}\label{covariance2}
\xi(gm)=\chi_l^{-1}(m) \xi(g) \qquad
  \forall\, g\in G, \ \ m\in M,
\end{gather}
and having a compact support modulo $M$.

\begin{lemma} Let $K\in \Hi_\rho^c$, $l\in \Sigma$ and $g\in G$. The
function
\[
\Gamma/\Gamma\cap
M \rightarrow \dsc, \qquad \gamma (\Gamma\cap M)\mapsto
K(l)(g\gamma)
\]
 is well defined.
\end{lemma}

\begin{proof} Let $\gamma\in \Gamma$ and $\gamma'\in\Gamma\cap M$.
Since $K(l)\in \mathcal{C}_c(G/M, \chi_l)$ then $K(l)$ satisf\/ies the
covariance relation (\ref{covariance2}). Hence $
K(l)(g\gamma\gamma') = \chi_l^{-1}(\gamma') K(l)(g\gamma)$. On the
other hand, as $\Gamma\cap M= \exp(\dsz X_1)\cdots \exp(\dsz
X_{n-1})$ and $M$ is abelian, then $\log(\Gamma)\cap \frak m= \dsz
X_1\oplus \dots \oplus \dsz X_{n-1}$. It follows that
$\chi_l\vert_{\Gamma\cap M}=1$, in particular, $\chi_l(\gamma')=1$.
Then $ K(l)(g\gamma\gamma') = K(l)(g\gamma)$.
\end{proof}

Let
\begin{gather*}
 \Hi_\rho^0 =\{K\in \Hi_\rho^c:\ K(l)  \mbox{is a zero function everywhere, except
for f\/inite number of  } l\}.
\end{gather*}
The space $ \Hi_\rho^0$ is dense in $\Hi_\rho$.

\begin{lemma} For $K\in \Hi_\rho^0 $ and $g\in G$, the sum
\begin{gather}\label{inverse}
  \sum_{l\in \Sigma} \sum_{\gamma\in \Gamma/\Gamma\cap
M} K(l)(g\gamma)+ \sum _{l\in \nongen} K(l)(g), \qquad g\in
G,
\end{gather} has only finitely many nonzero terms.
\end{lemma}

\begin{proof} Let  $l\in \Sigma\cup \nongen$.
For $g\in G$, let
\[
S_g=\set{\gamma\in \Gamma:\ g\gamma\in \mathrm{supp}(K(l))}.
\]
Then there is an integer $n_{K(l)}$ independent of $g\in G$ such
that $S_g$ is the union of at most $n_{K(l)}$ cosets of~$\Gamma\cap M$ \cite[Lemma~3.2]{Cor4}. This observation shows that the sum
over $\Gamma/\Gamma\cap M$  in~(\ref{inverse}) has at most
$n_{K(l)}$ nonzero entries for each~$g\in G$. Consequently, the sum~(\ref{inverse}) is a sum of f\/inite terms.
\end{proof}

We def\/ine the  operator $\V$  on the dense subspace $ \Hi_\rho^0 $
by
\begin{gather*}
\V(K)(g) =  \sum_{l\in \Sigma} \sum_{\gamma\in \Gamma/\Gamma\cap
M} K(l)(g\gamma)+ \sum _{l\in \nongen} K(l)(g),\qquad g\in G.
\end{gather*}

\begin{proposition}
The operator $\V$ is the inverse of $\U$.
\end{proposition}
\begin{proof}
First, we observe that $\V(K)$ satisf\/ies  the covariance relation
(\ref{covarience}) in $\mathbf{L}^2(G/\Gamma)$. In fact, let $g\in
G$ and $\gamma_0\in \Gamma$, we have
\begin{gather*}
\V(K)(g\gamma_0)   =  \sum_{l\in \Sigma} \sum_{\gamma\in
\Gamma/\Gamma\cap M} K(l)(g\gamma_0\gamma)+ \sum _{l\in
\nongen} K(l)(g\gamma_0)\\
\phantom{\V(K)(g\gamma_0)   =}{}   (\mbox{in the f\/irst sum, we use the  change of variable  }
\gamma\mapsto \gamma_0^{-1}\gamma)\\
\phantom{\V(K)(g\gamma_0)}{}  =  \sum_{l\in \Sigma} \sum_{\gamma\in \Gamma/\Gamma\cap M}
K(l)(g\gamma)+ \sum _{l\in \nongen} K(l)(g\gamma_0).
\end{gather*}
 As $\gamma_0\in \Gamma$ and $l\in \nongen$ then
 $\chi_l(\gamma_0)=1$ (see property $(1)$ of Lemma~\ref{periode2}).
 Consequently we obtain
\[
\V(K)(g\gamma_0)=\V(K)(g).
\]
 Now, we calculate $\U\circ \V (K)$
for some $K\in  \Hi_\rho$.
  Let $K\in \Hi_\rho^0$. We distinguish the following two cases.

\textbf{Case 1.}  $l_0\in \Sigma$.
\begin{gather*}
\U\circ \V(K)(l_0)(e) = \U( \V(K))(l_0)(e)  = \int_{M/M\cap \Gamma} \V(K)(m) \chi_{l_0}(m) d\dot{m}\\
\qquad{} = \int_{M/M\cap \Gamma} \Bigg( \sum_{l\in \Sigma}
\sum_{\gamma\in \Gamma/\Gamma\cap M} K(l)(m\gamma)+ \sum
_{l\in \nongen} K(l) (m)\Bigg)\chi_{l_0}(m) d\dot{m}\\
\qquad{} = \int_{M/M\cap \Gamma} \Bigg( \sum_{l\in \Sigma}
\sum_{\gamma\in \Gamma/\Gamma\cap M} K(l)(\gamma \gamma
^{-1}m\gamma)+ \sum
_{l\in \nongen} K(l) (m)\Bigg)\chi_{l_0}(m) d\dot{m}\\
\qquad{} = \int_{M/M\cap \Gamma} \Bigg( \sum_{l\in \Sigma}
\sum_{\gamma\in \Gamma/\Gamma\cap M} K(l)(\gamma
)\chi_{l}^{-1}(\gamma ^{-1}m\gamma) + \sum
_{l\in \nongen} K(l) (e)\chi_l^{-1}(m)\Bigg)\chi_{l_0}(m) d\dot{m}\\
\qquad{} = \int_{M/M\cap \Gamma} \Bigg( \sum_{l\in \Sigma}
\sum_{\gamma\in \Gamma/\Gamma\cap M} K(l)(\gamma
)\overline{\chi_{\Ad^*\gamma l}(m)} + \sum
_{l\in \nongen} K(l) (e)\overline{\chi_l(m)}\Bigg)\chi_{l_0}(m) d\dot{m}\\
\qquad{} =   \sum_{l\in \Sigma} \sum_{\gamma\in \Gamma/\Gamma\cap M}
K(l)(\gamma )\int_{M/M\cap \Gamma} \overline{\chi_{\Ad^*\gamma
l}(m)}\chi_{l_0}(m) d\dot{m} \\
\qquad\quad{} + \sum _{l\in \nongen} K(l) (e)
\int_{M/M\cap \Gamma} \overline{\chi_l(m)}\chi_{l_0}(m) d\dot{m}.
\end{gather*}
The interchange of integration and summation is justif\/ied by the
fact that $K\in \Hi^0_\rho$. On the other hand, the integral
\[
\int_{M/M\cap \Gamma}
\overline{\chi_{\Ad^*\gamma l}(m)} \chi_{l_0}(m) d\dot{m} = \left\{
  \begin{array}{ll}
    1, & \hbox{if} \ \Ad^*\gamma l = l_0, \\
    0, & \hbox{otherwise.}
  \end{array}
\right.
\]
Since $l$ and $l_0$ belong to $\Sigma$, then
\[
\Ad^*\gamma l= l_0 \ \Leftrightarrow \ l = l_{0}\qquad \mbox{and} \qquad
\gamma(\Gamma\cap M)= \Gamma\cap M.
\] It follows that
\[
\sum_{l\in \Sigma} \sum_{\gamma\in \Gamma/\Gamma\cap M}
K(l)(\gamma )\int_{M/M\cap \Gamma} \overline{\chi_{\Ad^*\gamma
l}(m)}\chi_{l_0}(m) d\dot{m}=K(l_0)(e).
\] On the other hand, since
$l_0\not\in \nongen$ then
\[
\int_{M/M\cap \Gamma}
\overline{\chi_l(m)}\chi_{l_0}(m) d\dot{m}=0,
\] for every $l\in
\nongen$ and hence
\[
\sum _{l\in \nongen} K(l) (e)
\int_{M/M\cap \Gamma} \overline{\chi_l(m)}\chi_{l_0}(m)
d\dot{m}=0.
\] Finally, we obtain
\[
\U\circ \V(K)(l_0)(e)= K(l_0)(e).
\]

 \textbf{Case 2.}   $l_0\in \nongen$.
\begin{gather*}
\U\circ \V(K)(l_0)(e)  = \U( \V(K))(l_0)(e) = \int_{G/ \Gamma} \V(K)(g) \chi_{l_0}(g) d\dot{g}\\
\phantom{\U\circ \V(K)(l_0)(e)}{} = \int_{G/ \Gamma}\Bigg( \sum_{l\in \Sigma} \sum_{\gamma\in
\Gamma/\Gamma\cap M} K(l)(g\gamma)+ \sum _{l\in \nongen} K(l)
(g)\Bigg) \chi_{l_0}(g) d\dot{g}.
\end{gather*}

First remark that:
\begin{gather*}
\int_{G/ \Gamma}   \sum_{l\in \nongen} K(l) (g) \chi_{l_0}(g)
d\dot{g}  = \int_{G/ \Gamma} \sum_{l\in \nongen}
K(l) (e) \overline{\chi_{l}(g)} \chi_{l_0}(g )d\dot{g}\\
\hphantom{\int_{G/ \Gamma}   \sum_{l\in \nongen} K(l) (g) \chi_{l_0}(g) d\dot{g} }{}
 =  \sum_{l\in \nongen} K(l) (e) \int_{G/ \Gamma}
\overline{\chi_{l}(g)} \chi_{l_0}(g)
 d\dot{g}
  = K(l_0)(e).
\end{gather*}

Next
\begin{gather*}
 \int_{G/ \Gamma}  \sum_{l\in \Sigma} \sum_{\gamma\in
\Gamma/ \Gamma\cap M} K(l)(g \gamma)\chi_{l_0}(g) d\dot{g}
 =
\int_{\tore^n}  \sum_{l\in \Sigma} \sum_{\gamma\in \Gamma/
\Gamma\cap M} K(l)(\delta(t, x)\gamma) \chi_{l_0}(\delta(t, x))
dx dt \\
 \qquad {} =  \int_{\tore^n}  \sum_{l\in \Sigma} \sum_{\gamma\in \Gamma/
\Gamma\cap M} K(l)(\delta(t, 0)\gamma \gamma^{-1}
\delta(0,x)\gamma)\chi_{l_0}(\delta(t, 0))
 \chi_{l_0}(\delta(0, x)) dx dt\\
\qquad {} =  \int_\tore \int_{\tore^{n-1}}  \sum_{l\in \Sigma}
\sum_{\gamma\in \Gamma/\Gamma\cap M} K(l)(\delta(t, 0)\gamma)
\chi_l^{-1}(\gamma^{-1} \delta(0,x)\gamma) \chi_{l_0}(\delta(t,
0))\chi_{l_0}(\delta(0, x)) dx dt\\
\qquad{} =  \int_\tore  \sum_{l\in \Sigma} \sum_{\gamma\in
\Gamma/\Gamma\cap M}
K(l)(\delta(t,0)\gamma)\int_{\tore^{n-1}}\overline{\chi_{\Ad^*\gamma
l}(\delta(0, x))} \chi_{l_0}(\delta(0, x)) dx \chi_{l_0}(\delta(t,
0))dt =0.
\end{gather*}
Then $\U\circ\V (K)(l_0)(e)\! =\! K(l_0)(e)$.
We conclude that for every $K$ in $\Hi_\rho^0$, we have $\U\circ\V
(K)\! =\! K $.

Next, for $\xi\in \mathbf{L}^2(G/\Gamma)$ such that $\U(\xi)\in
\mathscr{H}_\rho^0$ we compute
\begin{gather*}
\V\circ \U(\xi)(e)  = \V(\U(\xi))(e)
 =  \sum_{l\in \Sigma} \sum_{\gamma\in \Gamma/\Gamma\cap M}
\U(\xi)(l)(\gamma)+ \sum
_{l\in \nongen}\U(\xi)(l)(e)\\
 \phantom{\V\circ \U(\xi)(e)}{}
 =  \sum_{l\in \Sigma} \sum_{\gamma\in \Gamma/\Gamma\cap M}
\int_{M/M\cap \Gamma}\xi(\gamma m) \chi_l(m) d\dot{m} + \sum
_{l\in \nongen}\U(\xi)(l)(e)\\
\phantom{\V\circ \U(\xi)(e)}{}
=  \sum_{l\in \Sigma} \sum_{\gamma\in \Gamma/\Gamma\cap M}
\int_{M/M\cap \Gamma} \xi(\gamma m \gamma^{-1} ) \chi_l(m) d\dot{m}
+ \sum
_{l\in \nongen}\U(\xi)(l)(e)\\
\phantom{\V\circ \U(\xi)(e)}{}
=  \sum_{l\in \Sigma} \sum_{\gamma\in \Gamma/\Gamma\cap M}
\int_{M/M\cap \Gamma} \xi( m  ) \chi_l( \gamma^{-1} m \gamma)
d\dot{m} + \sum
_{l\in \nongen}\U(\xi)(l)(e)\\
\phantom{\V\circ \U(\xi)(e)}{}
=  \sum_{l\in \Sigma} \sum_{\gamma\in \Gamma/\Gamma\cap M}
\int_{M/M\cap \Gamma} \xi( m  ) \chi_{Ad^* \gamma l}(  m ) d\dot{m}
+ \sum
_{l\in \nongen}\U(\xi)(l)(e)\\
\phantom{\V\circ \U(\xi)(e)}{}
=  \sum_{l\in \gen}  \int_{M/M\cap \Gamma} \xi( m  ) \chi_l( m )
d\dot{m} + \sum
_{l\in \nongen}\U(\xi)(l)(e)\\
\phantom{\V\circ \U(\xi)(e)}{}
=  \sum_{l\in \gen} \int_{\tore^{n-1}} \xi(\delta(0, x))
\chi_l(\delta(0, x))dx  + \sum _{l\in \nongen}\U(\xi)(l)(e).
\end{gather*}

On the other hand
\begin{gather*}
\sum _{l\in \nongen}\U(\xi)(l)(e) = \sum _{l\in \nongen}
\int_{G/\Gamma} \xi(g) \chi_l(g) d\dot{g}
 = \sum_{l\in \nongen}\int_{\tore^{n}}
    \xi(\delta(t, x))
 \chi_l(\delta(t, x)) dt dx \\
\hphantom{\sum _{l\in \nongen}\U(\xi)(l)(e)}{}
 = \sum_{l_1\in \nongen_1} \sum_{l_2\in \nongen_2} \int_{\tore^{n}}
    \xi(\delta(t, x))
 \chi_{l_1}(\delta(0, x)) \chi_{l_2}(\delta(t, 0) ) dt dx \\
\hphantom{\sum _{l\in \nongen}\U(\xi)(l)(e)}{}
(\text{where } \nongen_1 = \mathds{Z} X_{p+1}^*+\cdots+\mathds{Z}X_{n-1}^*,
     \nongen_2 =\mathds{Z}X_{n}^* \text{ and } l = l_1+l_2
      \in \nongen_1\oplus \nongen_2) \\
\hphantom{\sum _{l\in \nongen}\U(\xi)(l)(e)}{}
= \sum_{l_1\in \nongen_1} \int_{\tore^{n-1}}
    \xi(\delta(0, x))
\chi_{l_1}(\delta(0, x))   dx \\
%\hphantom{\sum _{l\in \nongen}\U(\xi)(l)(e)}{} (\text{by the Fourier inversion formula to the variable }t)\\
\hphantom{\sum _{l\in \nongen}\U(\xi)(l)(e)}{}
\overset{\mbox{\scriptsize $\begin{array}{@{}c@{}}\text{by the Fourier inversion} \\ \text{formula to the variable  $t$}\end{array}$}}{=} \sum_{l\in \nongen_1} \int_{\tore^{n-1}}
    \xi(\delta(0, x))
 \chi_{l}(\delta(0, x))  dx.
\end{gather*}

Then
\begin{gather*}
\V\circ \U(\xi)(e)  =  \sum_{l\in \gen} \int_{\tore^{n-1}}
\xi(\delta(0, x)) \chi_l( \delta(0, x)) dx  + \sum_{l\in \nongen_1}
\int_{\tore^{n-1}}
    \xi(\delta(0, x)) \chi_{l}(\delta(0, x))  dx\\
\hphantom{\V\circ \U(\xi)(e)}{}
 \overset{\text{by (\ref{disjoint-reunion})}}{=}
   \sum_{l\in \frak m_{\mathbb{Z}}^*} \int_{\tore^{n-1}}
\xi(\delta(0, x))
\chi_l( \delta(0, x)) dx
\overset{\mbox{\scriptsize $\begin{array}{@{}c@{}}\text{by the Fourier inversion} \\ \text{formula to the variable  $x$}\end{array}$}}{=} \xi(e). \tag*{\qed}
\end{gather*}\renewcommand{\qed}{}
\end{proof}

Finally, we obtain that
 $\U$ is well def\/ined and isometric and
has dense range. It therefore extends uniquely into an isometry of
$\Hi_\tau$ onto $\Hi_\rho$.  As  immediate consequences of the last
theorem are the following

\begin{corollary} We have the following decomposition
\begin{gather*}%\label{decomposition-special}
\tau \simeq \rho =   \bigoplus_{l \in \Sigma}\ind_M^G\chi_l +
\bigoplus_{l\in \nongen} \chi_l.
\end{gather*}
\end{corollary}

For the next corollary, we write $\# E$ to denote the cardinal
number of a set $E$.
\begin{corollary}
We keep the same notations as Theorem  {\rm  \ref{des-operator-special}}.
The multiplicity function $\mult(\pi_l)$ is given by $\mult(\pi_l) =
1 $ if $l \in \nongen$ and  $\mult(\pi_l) = \#[\Omega(\pi_l) \cap
\Sigma] $ if $l \in \Sigma$.
\end{corollary}

\begin{proof} Let $l\in A$, where $A=\Sigma$ or $\nongen$. The
multiplicity $\mult(\pi_l)$ is equal to the number of $l'\in A$ such
that $\pi_l\simeq\pi_{l'}$. Thus, by the Kirillov theory, we have
\[
\mult(\pi_l)=\#\{l'\in A:\ l'\in \Omega(\pi_l)\}= \#[\Omega(\pi_l) \cap
A].
\]
 Now, if $l\in \nongen$, the coadjoint orbit $\Omega(\pi_l)$
has only one element and therefore $\mult(\pi_l)=1$.
\end{proof}

\pdfbookmark[1]{References}{ref}
\LastPageEnding

\end{document}